\documentclass{tran-l}
\usepackage{graphicx}
\usepackage{float}

\theoremstyle{definition}

\theoremstyle{remark}

\numberwithin{equation}{subsection}
\begin{document}

\title {Solve thermal explosion model by central difference and newton iteration method}

\author{
  Wang, Xijian
  \and
  Zeng, Tonghua
  }

\address{School of Mathematics \& Computational Science, Wuyi University, Guangdong Province, China}
\email{wangxj1980426@163.com, \ hnsdzth@163.com}

\thanks{This work was completed with the support of Foundation of Guangdong Educational Committee (2014KQNCX161, 2014KQNCX162) and Guangdong Province Comprehensive reform pilot major-Information and Computational Science, Wuyi University.}

%\thanks{The author was also supported in part by the Research Council of Slovenia.}

%\subjclass{Primary 47A15; Secondary 46A32, 47D20}

%\date{February 15, 1995 and, in revised form, July 6, 1995.}

%\dedicatory{}

%\commby{Daniel J. Rudolph}

% -----------------------------------------------------------------------------

\begin{abstract}
  First we present the general equation form of a thermal explosion in a vessel with boundary values, later use central difference method and Newton iteration method to solve the relevant partial differential equations in one-dimensional and two-dimensional forms, finally we use numerical experiments to verify the order of convergence of the central difference scheme and provide the experiment results.
\end{abstract}

% -----------------------------------------------------------------------------
\maketitle
% -----------------------------------------------------------------------------

\section*{Introduction}\label{intro}

N.N.Semenov ~\cite {semenov1959some}, Ya.B.Zeldovich ~\cite{zeldovich1985mathematical}, D.A.Frank-Kamene{\"e}t{\`\i}ski{\ae}i ~\cite{frank1969diffusion}, O.M.Todes and P.V.Melent'ev ~\cite{todes1939theory}, A.G.Merzhanov and F.I.Dubovitsky ~\cite{merzhanov1966modern}, B.Gray ~\cite{gray1973critical} had done a lot of research work in the field of thermal explosion theory.

\smallskip

An example of a boundary value problem(BVP) describing the onset of a thermal explosion in a vessel $\Omega$ is the following ~\cite{zeldovich1985mathematical}:
\begin{equation}\label{diffusion}
 \left\{ \begin{array}{l}
  a{\nabla ^2}T + Qk(T) = 0, x\in\Omega;\\
  T(x) = {T_0}, x\in \partial\Omega\\
 \end{array} \right.
\end{equation}
where $T(x)$ is the temperature in the vessel($[T]=K$). Equation (\ref{diffusion}) is a reaction-diffusion equation and determines the balance between conduction(thermal diffusion) and heat production by chemical reactions. Other variables/parameters in the BVP (\ref{diffusion}) are the ambient temperature $T_0([T_0]=K)$, the thermal diffusivity $a([a]=m^{2}/s)$, the heat release parameter $Q([Q]=K)$ and the reaction rate $k=k(T)([k]=s^{-1})$, given by the expression

  \begin{equation*}
    k(T)=Ae^{-T_{a}/T},
  \end{equation*}
where $A$ is the pre-exponential factor ($[A]=s^{-1}$) and $T_a$ the activation temperature $([T_a]=K)$. It is customary to rewrite equation \ref{diffusion} in terms of the dimensionless variables $u$ and $x^{*}$, defined by
\begin{equation*}
  u(x^{*}):=\frac{T_a}{T_0}\frac{T(x)-T_0}{T_o}, x^{*}:=\frac{x}{\ell},
\end{equation*}
where $\ell$ is a characteristic dimension of the vessel. To first order approximation, We can easily obtain
\begin{equation}\label{dimensionless}
  {\nabla ^2}u + q{e^u} = 0,\;q = \frac{{QA{\ell ^2}{T_a}{e^{ - \theta }}}}{{a{T_0}^2}},\;\theta  = \frac{{{T_a}}}{{{T_0}}}.
\end{equation}
A solution of equation (\ref{dimensionless}) exists if $q$ is small enough, i.e.,$q<q^{*}$ for some critical value $q^{*}$, meaning that chemical heat production can be entirely balanced by conduction. However, for increasing $q$ the solution of equation (\ref{dimensionless}) suddenly does not exist anymore. This phenomenon is referred to as a \emph{thermal explosion}. Since (\ref{dimensionless}) is a nonliear problem, ~\cite{mattheij2005partial} suggests we should use one of the following approaches: Newton iteration, Gauss-Jacobi iteration and transient methods. In this paper, we would first use central difference scheme to discretize the nonlinear equation (\ref{dimensionless}), later use Newton iteration method \cite{deuflhard2011newton,kelley2003solving,ortega2000iterative} to solve the nonlinear system of equations.

The paper is organized as the following: in section \ref{onedimension}, the one-dimensional case of equation (\ref{dimensionless}) is proposed and the solution is given analytically, later we use central difference method and Newton iteration method as the numerical algorithm to solve this equation, specially we use numerical experiments to obtain the threshold value of parameter $q$ and the convergence order of the algorithm by Richardson extrapolation method; in section \ref{twodimension}, the two-dimensional case of equation (\ref{dimensionless}) is proposed and again we use the same algorithm as in section \ref{onedimension} to solve the equation numerically with different values of parameter $q$, the numerical experiment results are listed as a table and several figures; finally the conclusion is given in section \ref{conclusion}.

\section{Solve One-dimensional thermal Explosion Model}\label{onedimension}

Consider the following one-dimensional BVP of (\ref{dimensionless})
\begin{equation}\label{onedim}
 \left\{ \begin{array}{l}
  u''+qe^{u}=0, x\in(0,1),\\
  u'(0)=0,u(1)=0.\\
 \end{array} \right.
\end{equation}

In this section, we will give the analytic of equation (\ref{onedim}), later use central difference method and Newton iteration method to compute the numerical solution of this equation, finally use Richardson extrapolation to compute the convergence order of the numerical scheme.

\subsection{Analytic solution}
The solution of equation (\ref{onedim}) can be found by subsequently multiplying equation with $u'$, integrate, isolate $u'$ and integrate again. This solution is given by
\begin{equation} \label{analytic}
  u(x)=ln(\frac{2\mu^2}{q})-2ln(cosh(\mu x)),
\end{equation}
where the parameter $\mu$ satisfies the relation
\begin{equation*}
  cosh\mu=\sqrt{2/q}\mu.
\end{equation*}

\subsection{Central difference and Newton iteration Scheme}
In order to compute the numerical solution of equation (\ref{onedim}), we introduce the grid
\begin{equation*}
  x_{j}=(j-1)\Delta x, j=1,2,\cdots, M,
  \end{equation*}
where $\Delta x=1/(M-1)$ is the grid size. We first use standard central differences for all derivatives in equation (\ref{onedim}). Let $u_{j}$ denote the numerical approximation of $u(x_{j})$ and let $\textbf{u}:=(u_{1}, u_{2}, \cdots, u_{M-1})$ be the vector of unknowns. Specially, we introduce the virtual point $x_{0}$ and apply the central difference to discretise the Robin boundary condition at the left boundary as
\begin{equation*}
  \frac{u_{2}-u_{0}}{2\Delta x}=u'(0)=0
\end{equation*}
Thus we have
\begin{equation*}
  u_{0}=u_{2}.
\end{equation*}
The resulting nonlinear system of algebraic equations can be written in the form
\begin{equation}\label{nonlinear1}
 \textbf{ N(u)}:=\textbf{Au}+\textbf{f(u)}=\textbf{0},
\end{equation}
where
\begin{equation*}
  \textbf{f(u)}=q\cdot diag(e^{u_{1}},e^{u_{2}},\cdots,e^{u_{M-1}}), A=\frac{1}{\Delta x^{2}}\begin{pmatrix}
         -2 & 2 & 0 & \cdots & 0 \\
         1 & -2 & 1 & \ddots & \vdots \\
         0 & \ddots & \ddots & \ddots & 0 \\
         \vdots & \ddots & \ddots & \ddots & 1 \\
         0 & \cdots & 0 & 1 & -2 \\
       \end{pmatrix}.
\end{equation*}
Now we would use Newton iteration scheme to solve equation (\ref{nonlinear1}). First give a suitable initial guess. Solve
\begin{equation*}
  \left\{ \begin{array}{l}
  u''+q=0, x\in(0,1),\\
  u'(0)=0,u(1)=0,\\
 \end{array} \right.
\end{equation*}
we obtain
\begin{equation*}
  u^{0}=\frac{q}{2}(1-x^2).
\end{equation*}
The Newton iteration scheme to solve equation (\ref{nonlinear1}) is the following:

\begin{equation*}
  \begin{array}{l}
   given,\, u^{0},\\
solve\, J(u^{l})s^{l}=-N(u^{l}),\\
 u^{l+1}=u^{l}+s^{l},
   \end{array}
\end{equation*}
with
\begin{equation*}
  u^{0}=\frac{q}{2}(1-x^2), J(u)=\frac{\partial N(u)}{\partial u}=(\frac{\partial N_{i}(u)}{\partial u_{j}}), s^{l}=-\frac{N(u^{l})}{J(u^{l})}.
\end{equation*}

\subsection{Numerical scheme to compute threshold value of Parameter $q$}
The numerical experiment result shows that the threshold value of parameter $q$ is 0.878. The numerical scheme to compute threshold value of Parameter $q$ is the following:\\
\textbf{Input:}grid size $\Delta x$,maxit(maximum iteration),tol(tolerance).\\
\textbf{Parameters:}$q=0.87:0.001:0.88$.\\
\textbf{Output:}threshold value of parameter $q^{*}$.\\
 $u^{0}: =\frac{q}{2}(1-x^2);$\\
 for $l=: 1$ to length($q$) do\\
 solve\, $J(u^{l})s^{l}=-N(u^{l})$,\\
 $u^{l+1}=u^{l}+s^{l}$;\\
 $F=A*u^{l+1}+f(u^{l+1})$;\\
 $nF=norm(F,2)$;\\
 conv=($nF<tol|(l==maxit)$);\\
 end\\
 if $(nF>tol)$\\
   \;\;\;\;  $q*=q(l-1)$;\\
   \;\;\;\; break\\
 else\\
 \;\;\;\; continue\\
 end\\
 $q^{*}$\\

\subsection{Convergence order} Due to the nonlinearity of the problem, convergence is difficult to prove. Instead, we verify the order of convergence by numerical experiments. Verify the order of convergence of the central difference scheme by computing numerical solutions of equation (\ref{nonlinear1}) for several values of the grid size $\Delta x$. Using Richardson Extrapolation method, the order of convergence can be computed by
\begin{equation*}
  p=\frac{ln(\frac{\|u_{2}-u_{e}\|_{\infty}}{\|u_{2}-u_{e}\|_{\infty}})}{ln2},
\end{equation*}
where $u_e$ represents the analytic solution from equation (\ref{analytic}), $u_1,u_2$ represents numerical solutions from different grid sizes. Numerical experiments show that the order of convergence  is 2.

\section{Solve One-dimensional thermal Explosion Model}\label{twodimension}
We consider the following two-dimensional BVP:
\begin{equation}\label{twodim}
  \left\{ \begin{array}{l}
  \frac{\partial ^{2}u}{\partial{x}^{2}}+ \frac{\partial ^{2}u}{\partial{y}^{2}}+qe^{u}=0, x\in (0,\ell),y\in(0,1),\\
  u(0,y)=0, u(\ell, y)=g(y), y\in(0,1),\\
  \frac{\partial u}{\partial y}(x,0)= \frac{\partial u}{\partial y}(x,1)=0, x\in(0,\ell),\\
 \end{array} \right.
\end{equation}
representing a rectangular vessel of aspect ratio $\ell$ with the horizontal walls isolated. For the boundary function $g(y)$ we take
\begin{equation*}
  g(y)=\left\{ \begin{array}{l}
 0,\,\, if \,\,y\in [0,\frac{1}{2})\\
 1,\,\, if \,\,y\in [\frac{1}{2},1].\\

 \end{array} \right.
\end{equation*}

\subsection{Numerical scheme}
For the numerical scheme of equation (\ref{twodim}) we use the following grid
\begin{equation*}
  (x_{j},y_{k})=((j-1)\Delta x, (k-1)\Delta y), j=1,2,\cdots, M, k=1,2,\cdots,N,
  \end{equation*}
  with $\Delta x=\ell/(M-1)$ and $\Delta y=1/(N-1)$ the grid sizes in $x-$ and $y-$ direction, respectively. We denote the numerical approximation of $u(x_{j}, y_{k})$ by $u_{j,k}$.
Similar to one-dimensional case, we use standard central differences for all derivatives in equation (\ref{twodim}). The resulting nonlinear system of algebraic equations can be written in the form
\begin{equation}\label{nonlinear2}
 N(U)=AU+bb+qe^U=0,
\end{equation}
where
\begin{equation*}
  A_{N\times N}=\begin{pmatrix}
         B & \frac{2}{\Delta y^2}I & O & \cdots & O \\
         \frac{1}{\Delta y^2}I & B & \frac{1}{\Delta y^2}I & \ddots & \vdots \\
         O & \ddots & \ddots & \ddots & O \\
         \vdots & \ddots & \frac{1}{\Delta y^2}I & B & \frac{1}{\Delta y^2}I \\
         O & O & O & \frac{2}{\Delta y^2}I & B \\
       \end{pmatrix},
\end{equation*}
\begin{equation*}
   B_{(M-2)\times(M-2)}=diag(-2(\frac{1}{\Delta x^2}+\frac{1}{\Delta y^2}),\cdots,-2(\frac{1}{\Delta x^2}+\frac{1}{\Delta y^2}),0),
\end{equation*}
\begin{equation*}
  bb=\begin{pmatrix}
       b_1 \\
       b_2 \\
       \vdots \\
       b_N \\
     \end{pmatrix}, U=\begin{pmatrix}
       u_1 \\
       u_2 \\
       \vdots \\
       u_N \\
     \end{pmatrix},  b_{i}=\begin{pmatrix}
       0 \\
       0 \\
       \vdots \\
       \frac{1}{\Delta y^2}g(y_i) \\
     \end{pmatrix}, u_i=\begin{pmatrix}
       u_{2i} \\
       u_{3i} \\
       \vdots \\
       u_{(M-1)i} \\
     \end{pmatrix}, i=1,2,\cdots,N.
\end{equation*}
Now we would use Newton iteration scheme to solve equation (\ref{nonlinear2}). First give a suitable initial guess. Solve
\begin{equation*}
  \left\{ \begin{array}{l}
  AU+bb+q=0,\\
  e^{U}=1,\\
 \end{array} \right.
\end{equation*}
we obtain
\begin{equation*}
  U^{0}=-A^{-1}(bb+q\cdot \overrightarrow{e}).
\end{equation*}
The Newton iteration scheme to solve equation (\ref{nonlinear2}) is the following:

\begin{equation*}
  \begin{array}{l}
   given,\,  U^{0},\\
solve\, J(U^{l})s^{l}=-N(U^{l}),\\
 U^{l+1}=U^{l}+s^{l},
   \end{array}
\end{equation*}
with
\begin{equation*}
  U^{0}=-A^{-1}(bb+q\cdot \overrightarrow{e}),\, J(U)=A+q\cdot diag(e^{U}), \,s^{l}=-\frac{N(U^{l})}{J(U^{l})} .
\end{equation*}
\subsection{Numerical results}
Similar to numerical scheme for one-dimensional case, we use central difference method and Newton iteration method to solve equation (\ref{twodim}). Let parameter $q$ be 0.8, 0.5, 0.3 respectively, given $\ell=1$, grid sizes $\Delta x=\Delta y=0.1$, we could obtain the number of iterations, solutions and the residues as the following table and figures:\\
\begin{center}
  \begin{tabular}{|c|c|c|c|c|c|}
  \hline
  % after \\: \hline or \cline{col1-col2} \cline{col3-col4} ...
  $q$ & $\ell$ & $\Delta x$& $\Delta y$ & number of iterations & residue\\
  \hline
  0.8 & 1 & 0.1 & 0.1 & 3 & 1.37145172e-013 \\
  \hline
  0.5 & 1 & 0.1 & 0.1 & 2 & 1.25775777e-009 \\
  \hline
  0.3 & 1 & 0.1 & 0.1 & 2 & 2.03563927e-011 \\
    \hline
\end{tabular}
\end{center}

\begin{figure}[H]
  \centering
  \begin{minipage}[t]{.35\linewidth}
   \includegraphics[width=2in]{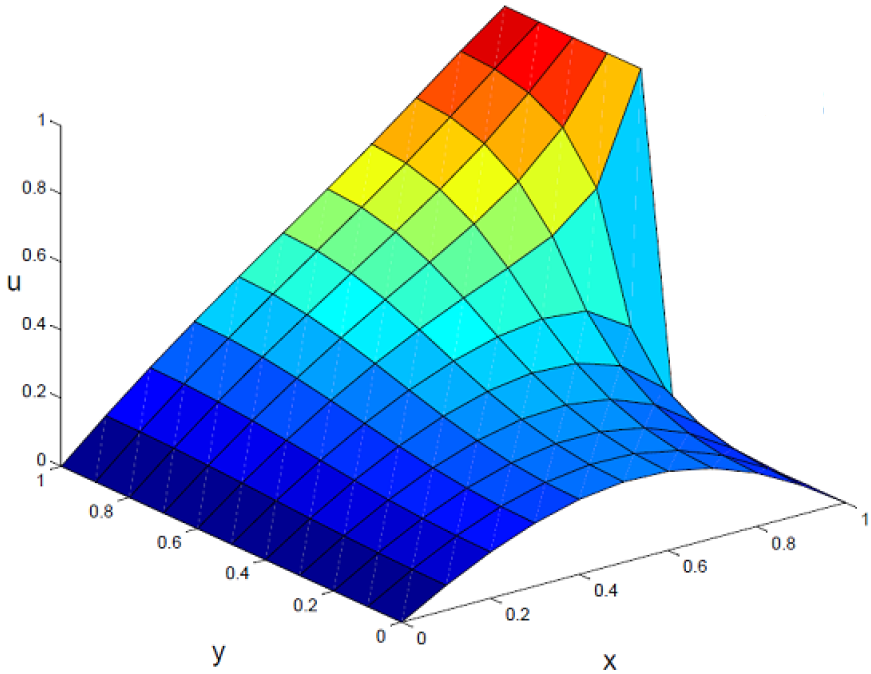}
  \caption{$q=0.8$}
  \end{minipage}
  \begin{minipage}[t]{.35\linewidth}
   \includegraphics[width=2in]{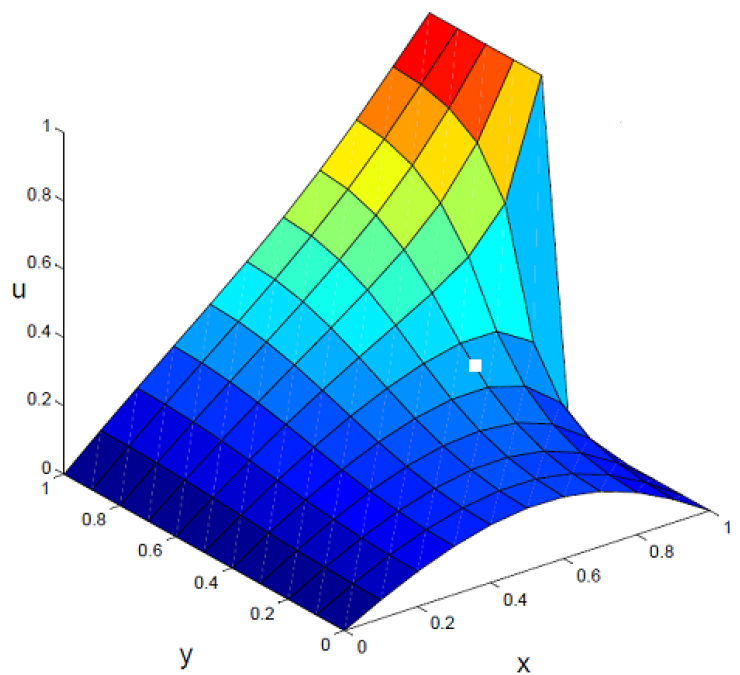}
  \caption{$q=0.5$}
  \end{minipage}
  \begin{minipage}[t]{.35\linewidth}
    \includegraphics[width=2in]{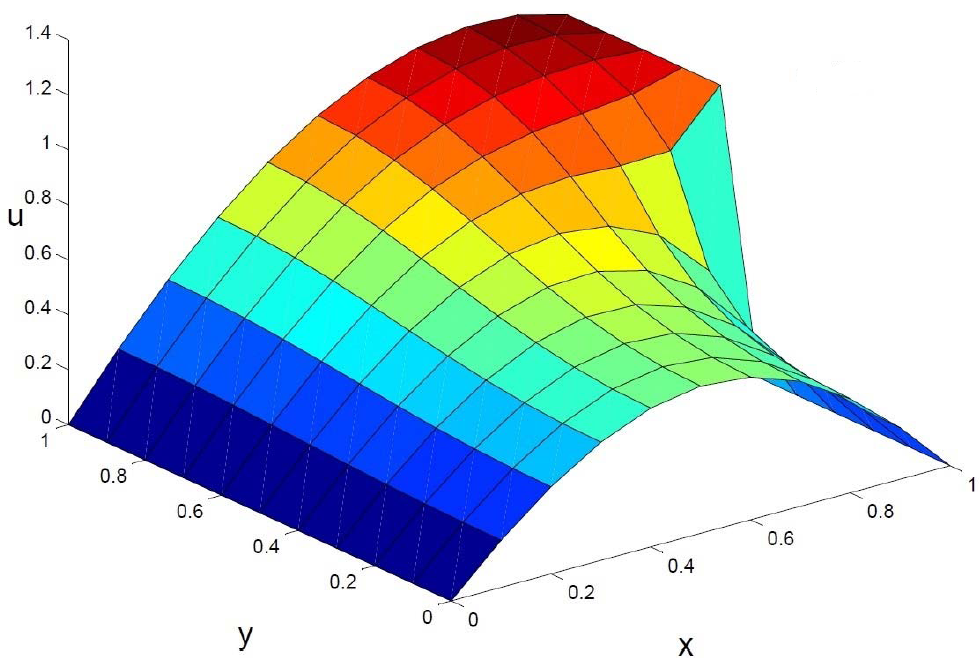}
  \caption{$q=0.3$}
  \end{minipage}
\end{figure}

\section{Conclusion}\label{conclusion}
We have considered a nonlinear boundary value problem describing the onset of a thermal explosion in a vessel in one-dimensional and two-dimensional forms. Numerical experiments show that Central different method combining Newton iteration method is suitable to solve the thermal explosion model (\ref{dimensionless}) if we find a suitable initial guess. In one-dimensional case (\ref{onedim}), Richardson extrapolation shows that the convergence order of our algorithm is 2; in two-dimensional case (\ref{twodim}), numerical experiments show that the number of iterations are less or equal to 3 and the residues are small.
\subsection*{Acknowledgment}
The results of this paper were obtained during our visiting study in Eindhoven University of Technology. we would like to express deep gratitude to our supervisor dr.ir. J.H.M. (Jan) ten Thije Boonkkamp and friend Dr Q. Hou whose guidance and support were crucial for the successful completion of this paper. 

% -----------------------------------------------------------------------------

\bibliographystyle{unsrt}

\begin{thebibliography}{10}

\bibitem{semenov1959some}
Nikola{\u\i}~Nikolaevich Semenov.
\newblock {\em Some problems of chemical kinetics and reactivity}, volume~2.
\newblock Elsevier, 1959.

\bibitem{zeldovich1985mathematical}
IA~Zeldovich, G~Io Barenblatt, VB~Librovich, and GM~Makhviladze.
\newblock {\em Mathematical theory of combustion and explosions}.
\newblock Consultants Bureau, New York, NY, 1985.

\bibitem{frank1969diffusion}
DA~Frank-Kamene{\"e}t{\`\i}ski{\ae}i.
\newblock {\em Diffusion and heat transfer in chemical kinetics}, volume~2.
\newblock Plenum Press (New York), 1969.

\bibitem{todes1939theory}
OM~Todes and PV~Melent'ev.
\newblock Theory of thermal explosion. 2. thermal explosion for unimoleclar
  reactions.
\newblock {\em Zh. Fiz. Khim}, 13:1594--1609, 1939.

\bibitem{merzhanov1966modern}
AG~Merzhanov and FI~Dubovitsky.
\newblock The modern state of the theory of thermal explosion.
\newblock {\em Uspekhi Khimii}, 35:4, 1966.

\bibitem{gray1973critical}
BF~Gray.
\newblock Critical behaviour in chemically reacting systems: Ii¡ªan exactly
  soluble model.
\newblock {\em Combustion and Flame}, 20(3):317--325, 1973.

\bibitem{mattheij2005partial}
Robert~MM Mattheij, Sjoerd~W Rienstra, and Jan~HM ten Thije~Boonkkamp.
\newblock {\em Partial differential equations: modeling, analysis,
  computation}.
\newblock Siam, 2005.

\bibitem{deuflhard2011newton}
Peter Deuflhard.
\newblock {\em Newton methods for nonlinear problems: affine invariance and
  adaptive algorithms}, volume~35.
\newblock Springer, 2011.

\bibitem{kelley2003solving}
Carl~T Kelley.
\newblock {\em Solving nonlinear equations with Newton's method}, volume~1.
\newblock Siam, 2003.

\bibitem{ortega2000iterative}
James~M Ortega and Werner~C Rheinboldt.
\newblock {\em Iterative solution of nonlinear equations in several variables},
  volume~30.
\newblock Siam, 2000.

\end{thebibliography}

\end{document}